\documentclass[a4paper, 11pt]{amsart}

\usepackage[mathscr]{eucal}
\usepackage{amsmath}
\usepackage{amssymb}
\usepackage{amsfonts}
\usepackage{amsthm}

\usepackage{cases}

\usepackage[mathscr]{eucal}
\usepackage{amsmath}
\usepackage{amssymb}
\usepackage{amsfonts}
\usepackage{amsthm}
\usepackage[dvips]{graphics, color}
\theoremstyle{plain}
\newtheorem{theorem}{Theorem}[section]
\newtheorem{lemma}{Lemma}[section]
\newtheorem{remark}{Remark}[section]

\renewcommand{\labelenumi}{(\theenumi)}
\numberwithin{equation}{section}

\def\<{\left<} \def\>{\right>}

\def\bea{\begin{eqnarray} }
\def\eea{\end{eqnarray} }
\def\be{\begin{equation} }
\def\ee{\end{equation} }
\def\qed{\ifhmode\unskip\nobreak\fi\ifmmode\ifinner\else\hskip5pt
\fi\fi\hbox{\hskip5 pt \vrule width4 pt height6 pt depth1.5 pt \hskip1pt }}

\begin{document}
\title[]{Real hypersurfaces in complex space forms attaining equality in an inequality
involving a contact $\delta$-invariant}
\author[]{Toru Sasahara}
\address{Division of Mathematics, 
Hachinohe Institute of Technology, 
Hachinohe, Aomori 031-8501, Japan}
\email{sasahara@hi-tech.ac.jp}


\begin{abstract} 
We investigate real hypersurfaces in complex space forms  attaining equality in an inequality
involving the contact $\delta$-invariant $\delta^c(2)$ introduced by
Chen and Mihai in \cite{chen2}.
\end{abstract}



\maketitle

 \section{Introduction}
 Let  $\tilde M^n(4c)$ denote an $n$-dimensional complex space form
 of constant holomorphic sectional curvature $4c$ $(\ne 0)$, that is, the complex projective $n$-space
 $\mathbb{C}P^n(4c)$ or the complex hyperbolic space $\mathbb{C}H^n(4c)$, according as
 $c>0$ or $c<0$.
 We denote by $J$  the almost complex structure on $\tilde M^n(4c)$.
 
 Let $M$ be a real hypersurface of  $\tilde M^n(4c)$.
 For a unit normal vector field $N$ of  $M$ in $\tilde M^n(4c)$,
the structure vector field on $M$ is defined by $\xi=-JN$.
We define  a $1$-form $\eta$ and a $(1, 1)$-tensor $\phi$
 by $\eta(X)=g(\xi, X)$ and $\phi X=JX-\eta(X)N$ for each vector field $X$ tangent to $M$,
  where $g$ is the induced Riemannian metric.
 Then, a quadruplet $(\phi, \xi, \eta, g)$ defines 
 an almost contact metric structure on
 $M$.
If $\xi$ is a principal curvature vector everywhere on $M$, 
then $M$ is called a {\it Hopf hypersurface}.

Let $\mathcal{H}$ be the holomorphic distribution defined by $\mathcal{H}_p=\{X\in T_pM\ |\  \eta(X)=0\}$ for $p\in M$.
If $\mathcal{H}$ is integrable and each leaf of its maximal integral manifolds is 
locally congruent to
a totally geodesic complex hypersurface
$\tilde M^{n-1}(4c)$ in $\tilde M^n(4c)$, then $M$ is called a {\it ruled real hypersurface}, which is a
 typical example of a non-Hopf hypersurface.



For a Riemannian manifold $M$ endowed with  an almost contact metric structure 
$(\phi, \xi, \eta, g)$, Chen and Mihai \cite{chen2} defined
the contact $\delta$-invariant $\delta^c(2)$ of $M$  
by $\delta^c(2)(p)=\tau(p)-\inf\{K(\pi_{\xi})\ |\ \pi_{\xi}\  \text{is a plane containing $\xi$ in}\ T_pM\}$,
where $\tau$ is the scalar curvature of $M$ and $K(\pi_{\xi})$ is the sectional curvature of $\pi_{\xi}$. 

From Theorem 2 in \cite{chen1}, it immediately follows that
 a real hypersurface in $\tilde M^n(4c)$ satisfies
\be
   \delta^c(2)\leq \frac{(2n-1)^2(2n-3)}{4(n-1)}\|H\|^2+(2n^2-3)c \label{ineq}
   \ee 
   for the induced almost contact metric structure $(\phi, \xi, \eta, g)$, where
   $H$ denotes the mean curvature vector.
A real hypersruface is called {\it $\delta^c(2)$-ideal} if it attains
equality in (\ref{ineq}) 
at every point.  
 

 In this paper, we obtain the following two classification results for $\delta^c(2)$-ideal
 hypersurfaces in $\tilde M^n(4c)$.
 \begin{theorem}\label{Hopf}
 Let $M$ be a Hopf hypersurface in a complex space form $\tilde M^n(4c)$,
 where $c\in\{-1, 1\}$.
 Then  $M$ is $\delta^c(2)$-ideal
  if and only if $M$ is locally congruent to one of the following{\rm :}
 \begin{enumerate}
\renewcommand{\labelenumi}{(\roman{enumi})}
 \item a geodesic sphere of radius $\pi/4$ in  $\mathbb{C}P^n(4)$,
 
 \item a tube of radius  $r=\tan^{-1}((1+\sqrt{5}-\sqrt{2+2\sqrt{5}})/2)$
 around a complex quadric curve $Q_1$ in $\mathbb{C}P^2(4)$,

 \item a tube of radius 
  $(1/2)\log((1+\sqrt{5}\pm\sqrt{2+2\sqrt{5}})/2)$ around a totally real totally geodesic $\mathbb{R}H^2$
  in
  $\mathbb{C}H^2(-4)$. 
\end{enumerate} 
 
 \end{theorem}

\begin{theorem}\label{Non-Hopf}
Let $M$ be a non-Hopf real hypersurface  in a complex space form.
Assume that $M$ has constant mean curvature.  Then 
$M$ is $\delta^c(2)$-ideal if and only if 
$M$ is a minimal ruled real hypersurface.
\end{theorem}

\begin{remark}
{\rm Minimal ruled real hypersurfaces in non-flat complex space forms have
been classified in \cite{ada}.}
\end{remark}

\begin{remark}
{\rm 
Suppose that $\alpha(s)$, $\beta(s)$, $\gamma(s)$ and $\mu(s)$ satisfy the following system of ODE's and algebraic condition:
\be
\begin{split}
& \frac{d\alpha}{ds}=\beta(\alpha+\gamma-3\mu),\\
& \frac{d\beta}{ds}=\beta^2+\gamma^2+\mu(\alpha-2\gamma)+c,\\
& \frac{d\gamma}{ds}=\frac{(\gamma-\mu)(\gamma^2-\alpha\gamma-c)}{\beta}
+\beta(2\gamma+\mu),\\
& \alpha+\gamma=\mu, 
\end{split} \nonumber
\ee
where $\mu$ and $\beta$ are non-vanishing.
Then there exists  a non-minimal non-Hopf real hypersurface in $\tilde M^2(4c)$,
such that the components of the shape operator are given by (\ref{A}) in Section $2$, where
$d/ds=e_3$ (see Theorem 5 in \cite{ivey}). 
}
\end{remark}

\section{Preliminaries}

Let $M$ be a   real hypersurface in a complex space form $\tilde M^n(4c)$. 
Let  $(\phi, \xi, \eta, g)$ be an almost contact metric structure induced from the complex
structure of  $\tilde M^n(4c)$. 
Let us denote by $\nabla$ and
 $\tilde\nabla$ the Levi-Civita connections on $M$ and $\tilde M^n(4c)$, respectively. The
 Gauss and Weingarten formulas are,  respectively,  given by
\be
 \begin{split}
 \tilde \nabla_XY&= \nabla_XY+g(AX, Y)N, \label{gawe}\\
 \tilde\nabla_X N&= -AX
 \end{split}\nonumber
\ee
 for tangent vector fields $X$, $Y$ and a unit normal vector field $N$,
 where $A$ is the shape operator.
The mean curvature vector field $H$ is defined by 
$H=({\rm Tr}A/(2n-1))N.$
 The function ${\rm Tr}A/(2n-1)$ is called the  {\it mean curvature}.
 If it vanishes identically, then $M$ is called a {\it minimal hypersurface}.
 
By the Gauss and  Weingarten formulas, we have
\be
\nabla_{X}\xi=\phi AX. \label{PA}
\ee

 We denote by $R$ the Riemannian curvature tensor of $M$. Then,
 the equations of Gauss  and Codazzi are, respectively, given by
 \begin{gather}
R(X, Y)Z=c[g(Y, Z)X-g(X, Z)Y+g(\phi Y, Z)\phi X
 -g(\phi X, Z)\phi Y]  \label{ga}\\
  -2g(\phi X, Y)\phi Z
   +g(AY, Z)AX-g(AX, Z)AY,\nonumber\\
  (\nabla_XA)Y-(\nabla_YA)X=c[\eta(X)g(\phi Y, Z)-
  \eta(Y)g(\phi X, Z)-2g(\phi X, Y)\xi]. \label{co}
 \end{gather}

The following two lemmas are crucial.
 \begin{lemma}[\cite{chen1}]\label{lem1}
  Let $M$ be a real  hypersurface 
   in $\tilde M^n(4c)$. Then the equality sign in $(\ref{ineq})$ holds at a point $p\in M$ 
if and only if there exists an orthonormal basis $\{e_1, \dots, e_{2n-1}\}$ at $p$ 
 such that $e_1=\xi$, $e_{2i+1}=\phi e_{2i}$ $(1\leq i \leq n-1)$
 and the shape operator of $M$ in $\tilde M^n(4c)$ at $p$ satisfies
\be
A= \left(
    \begin{array}{ccccc}
      \alpha & \beta & 0 & \ldots & 0 \\
      \beta & \gamma & 0 & \ldots & 0 \\ 
      0 & 0 & \mu &\dots & 0 \\
      \vdots & \vdots & \vdots & \ddots & \vdots \\
      0 & 0 & 0 & \ldots & \mu
   \end{array}
  \right), \ \ \alpha+\gamma=\mu. \label{A} 
\ee
\end{lemma}

\begin{lemma}[\cite{ni}]\label{lem2}
Let $M$ a  real hypersurface $M$ in $\tilde M^n(4c)$ with $n\geq 2$.
We define differentiable functions $\alpha$, $\beta$ on $M$ by $\alpha=g(A\xi, \xi)$
and $\beta=\|A\xi-\alpha\xi\|$.
Then $M$ is ruled if and only if 
the following two conditions hold{\rm :} 




{\rm (1)} the set $M_1=\{p\in M \ | \ \beta(p)\ne 0\}$ is an open dense subset of $M${\rm ;}

{\rm (2)} there is a unit vector field $U$ on $M_1$, which is orthogonal to $\xi$ and satisfies
\be 
A\xi=\alpha\xi+\beta U, \ \ AU=\beta\xi, \ \ AX=0 \nonumber \label{ruled}
\ee
 for an arbitrary tangent  vector $X$ orthogonal to both $\xi$ and $U$.

\end{lemma}

\section {Proof of Theorem 1.1}
The proof of Theorem 8 in \cite{chen1} shows that
a $\delta^c(2)$-ideal Hopf hypersurface in $\mathbb{C}P^n(4)$
 is locally congruent to (i) or (ii).
 In order to investigate the case $c=-1$,  we will first describe
 some well-known fundamental results regarding Hopf hypersurfaces  
 in  $\mathbb{C}H^n(-4)$ (cf. \cite{ni}).
 \begin{theorem}\label{eigen}
 Let $M$ be a Hopf hypersurface in  $\mathbb{C}H^n(-4)$ 
 and $\nu$ 
 the principal curvature corresponding to the characteristic vector field $\xi$. 
 Then we have
 \begin{itemize}
 \item[(i)] $\nu$ is constant{\rm ;}
 \item[(ii)] If $X$ is a  tangent vector of $M$ orthogonal to $\xi$  such that 
  $AX=\lambda_1X$ and $A\phi X=\lambda_2\phi X$, then
 $2\lambda_1\lambda_2=(\lambda_1+\lambda_2)\nu-2$ holds.
 \end{itemize} 
 \end{theorem}

\begin{theorem}\label{hopf}
Let $M$ be a Hopf hypersurface with constant principal curvatures in $\mathbb{C}H^n(-4)$.
 Then $M$ is  locally congruent to one of the following: 
\begin{itemize}
\item[$(A_0)$] a horosphere, 
\item[$(A_{1,0})$] a geodesic sphere of radius $r$, where $0<r<\infty$, 
\item[$(A_{1,1})$] a tube of  radius $r$ around a totally geodesic $\mathbb{C}H^{n-1}(-4)$, 
where $0<r<\infty$, 
\item[$(A_2)$] a tube of radius $r$ around a totally geodesic $\mathbb{C}H^k(-4)$, where
 $1\leq k\leq n-2$ and $0<r<\infty$,
\item[$(B)$] a tube of radius $r$ around a totally real totally geodesic $\mathbb{R}H^n$, where $0<r<\infty$.
\end{itemize}
\end{theorem}

\begin{theorem}\label{A0}
The Type $(A_0)$ hypersurfaces in $\mathbb{C}H^n(-4)$
have two principal curvatures: 
$\nu=2$ of multiplicity $1$ and
$\lambda_1=1$ 
of multiplicity $2n-2$.
\end{theorem}

\begin{theorem}\label{A10}
The Type $(A_{1,0})$ hypersurfaces in $\mathbb{C}H^n(-4)$ have two distinct principal curvatures: 
 $\nu=2\coth(2r)$  of multiplicity $1$ and 
$\lambda_1=\coth(r)$
of multiplicity $2n-2$.
\end{theorem}

\begin{theorem}\label{A11}
The Type $(A_{1,1})$ hypersurfaces
in  $\mathbb{C}H^n(-4)$ have two distinct principal curvatures: 
 $\nu=2\coth(2r)$  of multiplicity $1$ and 
$\lambda_1=\tanh(r)$
of multiplicity $2n-2$. 
\end{theorem}

\begin{theorem}\label{A2}
The Type $(A_2)$ hypersurfaces
 in $\mathbb{C}H^n(-4)$ have three  distinct principal curvatures: 
$\nu=2\coth(2r)$ of multiplicity $1$,
$\lambda_1=\coth(r)$ of multiplicity $2n-2k-2$, and 
$\lambda_2=\tanh(r)$ of multiplicity $2k$.
\end{theorem}

\begin{theorem}\label{B}
The Type $(B)$ hypersurfaces
 in $\mathbb{C}H^n(-4)$ have three  principal curvatures: 
$\nu=2\tanh(2r)$ of multiplicity $1$,
$\lambda_1=\coth(r)$ of multiplicity $n-1$, and 
$\lambda_2=\tanh(r)$ of multiplicity $n-1$.
These principal curvatures are distinct unless
$r=(1/2)\log(2+\sqrt{3})$, in which case $\nu=\lambda_1=\sqrt{3}$ and 
$\lambda_2=1/\sqrt{3}$.
\end{theorem}

Let $M$ be a Hopf hypersurface in $\mathbb{C}H^n(-4)$ attaining 
equality in $(\ref{ineq})$ at every point.
By Lemma \ref{lem1},    
the  shape operator of $M$ takes the form (\ref{A}) with $\beta=0$ 
 at each point.  
 We put $\nu=\alpha$, $\lambda_1=\gamma$ and $\lambda_2=\mu$
 in Theorem \ref{eigen}. It follows from Theorem \ref{eigen} and 
 (\ref{A}) 
that $M$ is a Hopf hypersurface with constant principal curvatures.
Thus, we can apply Theorems \ref{hopf}-\ref{B}.  

Taking into account (\ref{A}) and the multiplicities of principal curvatures,
we see that $M$ is  a Type $B$ hypersurface in $\mathbb{C}H^2(-4)$.
  Since  $-1<\tanh(r)<1$,  the equation $2\tanh(2r)+\tanh(r)=\coth(r)$
  has no solutions, and the equation
 $2\tanh(2r)+\coth(r)=\tanh(r)$  has solutions 
  $(1/2)\log((1+\sqrt{5}\pm\sqrt{2+2\sqrt{5}})/2)$. 
 
 The converse is clear from Lemma \ref{lem1}.
  The proof of Theorem \ref{Hopf} is finished.

\section{Proof of Theorem 1.2}
 Let $M$ be a non-Hopf real hypersurface  in a complex space form $\tilde M^n(4c)$,
 where $c\in\{-1, 1\}$.
Assume that $M$ has constant mean curvature 
and attains equality in $(\ref{ineq})$ at every point.

   Let $\{e_1, \ldots, e_{2n-1}\}$ be a local orthonormal frame field described in Lemma \ref{lem1}. 
Then, (\ref{A}) can be rewritten as
\begin{align} 
& A\xi=(\mu-\gamma)\xi+\beta e_2,\ \  Ae_2=\gamma e_2+\beta\xi, 
\ \ Ae_i=\mu e_i,\label{A2}
 \end{align}
 where $i=3, \ldots, 2n-1$.
Since $\xi$ is not a principal vector everywhere, 
we have $\beta\ne 0$ on $M$.
The constancy of the mean curvature  implies that $\mu$ is constant.
By (\ref{PA}) and (\ref{A}), we get
\be
\nabla_{e_2}\xi=\gamma e_3, \ \ \nabla_{e_3}\xi=-\mu e_2, \ \ \nabla_{\xi}\xi=\beta e_3. \label{n1}
\ee
We put
$\kappa_1=g(\nabla_{e_2}e_2, e_3)$, $\kappa_2=g(\nabla_{e_3}e_2, e_3)$ and  $
\kappa_3=g(\nabla_{\xi}e_2, e_3)$.

We compare
 the coefficients with respect to $\{\xi, e_2, e_3\}$ on both sides of 
  the equation  (\ref{co}) of Codazzi for $X, Y\in\{\xi, e_2, e_3\}$.
 Then, taking into account  $g(\nabla{e_i}, e_j)=-g(\nabla{e_j}, e_i)$, 
 and using (\ref{A2}), (\ref{n1}), we get
  \begin{align}
   e_2\gamma &=-\xi\beta, \label{cd6}\\
   e_2\beta &=\xi\gamma, \label{cd7}\\
  \beta\kappa_1+(\mu-\gamma)\kappa_3 &=\beta^2+\gamma^2-c,\label{cd1}\\
e_3\gamma &=2\beta\mu+\beta\gamma-\beta k_3,\label{cd3}\\
e_3\beta &=\mu^2-2\mu\gamma-\kappa_3(\mu-\gamma)+\beta^2+c, \label{cd2}\\
   \kappa_2 &=0,\label{cd4}\\
e_3\gamma &=(\gamma-\mu)\kappa_1+2\beta\mu+\beta\gamma,\label{cd5}\\
 e_3\beta &=\kappa_1\beta+\mu^2-\gamma^2+2c. \label{cd9}
\end{align}

Eliminating $e_3\gamma$ from (\ref{cd3}) and (\ref{cd5}), we obtain
 \be
 (\mu-\gamma)\kappa_1-\beta\kappa_3=0.
 \label{eq1}
 \ee 
  Solving (\ref{cd1}) and (\ref{eq1}) for $\kappa_1$ and $\kappa_3$ yields
 \begin{align}
 \kappa_1 &=\frac{\beta(\beta^2+\gamma^2-c)}{(\mu-\gamma)^2+\beta^2}, \label{k1}\\
 \kappa_3 &=\frac{(\mu-\gamma)(\beta^2+\gamma^2-c)}{(\mu-\gamma)^2+\beta^2}.\label{k3}
 \end{align}
Substituting (\ref{k3}) and (\ref{k1}) into (\ref{cd2}) and (\ref{cd5}), respectively,  we obtain
 \begin{align}
 &e_3\beta=(\mu-2\gamma)\mu+\beta^2+c-\frac{(\mu-\gamma)^2(\beta^2+\gamma^2-c)}{(\mu-\gamma)^2+\beta^2}, \label{b}\\
& e_3\gamma=\beta(\gamma+2\mu)+\frac{(\gamma-\mu)\beta(\beta^2+\gamma^2-c)}{(\mu-\gamma)^2+\beta^2}.\label{r}
\end{align}
Note that substitution of (\ref{k1}) into (\ref{cd9}) gives (\ref{b}).

{\bf Case (I)}. \ $n=2$.

If $c=1$, then the statement of Theorem 1.2 can be proved in the same way as in the proof of Theorem 1.2 in \cite{sa}. Hence we consider only the case $c=-1$.

Using the equation (\ref{ga}) of Gauss for $g(R(e_2, e_3)e_3, e_2)$, and
taking into account $(\ref{cd4})$,
 we obtain 
 \be
 e_3\kappa_1-2\mu\gamma-\kappa_1^2-(\gamma+\mu)\kappa_3+4=0.\label{gauss}
 \ee
 Substituting (\ref{k1}) and (\ref{k3}) into (\ref{gauss}) gives
 \begin{align}
 &\{(3\beta^2+\gamma^2+1)((\mu-\gamma)^2+\beta^2)-2\beta^2(\beta^2+\gamma^2+1)\}e_3\beta  \label{g1}\\
 &+\{2\beta\gamma((\mu-\gamma)^2+\beta^2)+2(\mu-\gamma)\beta(\beta^2+\gamma^2+1)\}e_3\gamma\nonumber \\
 &-2\mu\gamma\{(\mu-\gamma)^2+\beta^2\}^2-\beta^2(\beta^2+\gamma^2+1)^2\nonumber\\
 &+(\gamma^2-\mu^2)(\beta^2+\gamma^2+1)\{(\mu-\gamma)^2+\beta^2\}+4\{(\mu-\gamma)^2+\beta^2\}^2=0.\nonumber
 \end{align}
 Substituting   (\ref{b}) and   (\ref{r}) into   (\ref{g1}), we have
 \be
 (\mu-\gamma)f(\beta, \gamma)=0, \label{f}\ee
where  $f(\beta, \gamma)$ is a polynomial given by 
\be
\begin{split}
f(\beta, \gamma):=& 2\mu\gamma^4-(4\mu^2+1)\gamma^3+(3\mu^2+4\beta^2+6)\mu\gamma^2-\{\mu^4+(4\beta^2+7)\mu^2+\beta^2-1\}\gamma\nonumber\\
&+(\beta^2+2)\mu^3+(2\beta^4+2\beta^2-1)\mu.\nonumber
\end{split}
\ee
 
 
 {\bf Case (I.1)}. \  $\mu-\gamma=0$.
 
   In this case,  by (\ref{eq1}) we get $\kappa_3=0$. 
 Therefore, by (\ref{cd3})  and the constancy of $\mu$, we find
  that $\mu=\gamma=0$.

{\bf Case (I.2)}. \ $f(\beta, \gamma)=0$.

We differentiate $f(\beta, \gamma)=0$ along $e_3$. Then,
 using (\ref{b}) and (\ref{r}),  we obtain
\be
\begin{split}
& 8\mu\gamma^6-(24\mu^2+4)\gamma^5+(30\mu^2+24\beta^2+15)\mu\gamma^4\\
 &-\{20\mu^4+(48\beta^2-3)\mu^2+8b^2-3\}\gamma^3\\
 &+\{7\mu^5+(36\beta^2-45)\mu^3+(24\beta^4+10\beta^2-2)\mu\}\gamma^2\\
&-\{\mu^6+(12\beta^2-44)\mu^4+(24\beta^4-19\beta^2+2)\mu^2+4\beta^4-3\beta^2-1\}\gamma\\
 &+(\beta^2-13)\mu^5+(6\beta^4-19\beta^2+1)\mu^3+(8\beta^6+8\gamma^6-5\beta^4-2\beta^2-1)\mu=0. \label{f2}
\end{split}
\ee
 The resultant of $f(\beta, \gamma)$ and the left-hand side of (\ref{f2}) with respect to $\gamma$ is given by
  \be
  202500(\mu^2+1)^4\beta^4\mu^6\{4\mu^2\beta^2+(\mu^2+1)^2\}^2.\nonumber
  \ee
Since $\beta\ne 0$, we have $\mu=0$.
 Therefore, equation $f(\beta, \gamma)=0$ can be simplified to
 \be 
 \gamma(\beta^2+\gamma^2-1)=0. \nonumber\label{eq2}
 \ee

 We will investigate the case $\beta^2+\gamma^2=1$.
 It follows from (\ref{k1}),  (\ref{k3}), (\ref{b}) and (\ref{r}) that 
 \begin{align} 
& \kappa_1=2\beta, \ \ \kappa_3=-2\gamma, \label{k13}\\
 & e_3\beta=3(\beta^2-1), \ \ e_3\gamma=3\beta\gamma. \label{br}
 \end{align}
 Differentiating $\beta^2+\gamma^2=1$ along $e_3$, we have
$\beta e_3\beta+\gamma e_3\gamma=0$. However,  substitution of (\ref {br}) into this equation
provides no new information. 

Differentiating $\beta^2+\gamma^2=1$ along $\xi$ gives
\be \beta(\xi\beta)+\gamma(\xi\gamma)=0. \label{eq5}\ee
Moreover, differentiating (\ref{eq5}) along $e_3$, we obtain
\be 
(e_3\beta)(\xi\beta)+\beta e_3(\xi\beta)+(e_3\gamma)(\xi\gamma) \label{eq6}
+\gamma e_3(\xi\gamma)=0.
\ee 
Using (\ref{n1}), (\ref{cd6}), (\ref{cd7}), (\ref{k13}) and (\ref{br}),
 we obtain
\begin{align}
e_3(\xi\beta) &=(\nabla_{e_3}\xi-\nabla_{\xi}e_3)\beta+\xi(e_3\beta) \nonumber \\
&=(\kappa_3e_2+\beta\xi)\beta+6\beta(\xi\beta),\nonumber\\
&=7\beta\xi\beta-2\gamma(\xi\gamma), \label{eq3}\\
e_3(\xi\gamma) &=(\nabla_{e_3}\xi-\nabla_{\xi}e_3)\gamma+\xi(e_3\gamma)  \nonumber \\
&=(\kappa_3e_2+\beta\xi)\gamma+3\gamma(\xi\beta)+3\beta(\xi\gamma)\nonumber \\
&=5\gamma(\xi\beta)+4\beta(\xi\gamma).\label{eq4}
\end{align}
Substitution of  (\ref{br}), (\ref{eq3}) and (\ref{eq4}) into (\ref{eq6}) implies
\be 
(10\beta^2+5\gamma^2-3)(\xi\beta)+5\beta\gamma(\xi\gamma)=0.\label{eq7}
\ee
Eliminating $\xi\beta$ from (\ref{eq5}) and (\ref{eq7}), and using $\beta^2+\gamma^2=1$,
we get
\be 
\gamma(\xi\gamma)=0.
\ee 
Thus, $\xi\gamma=0$. It follows from (\ref{cd6}), (\ref{cd7}) and (\ref{eq5}) that
$\xi\beta=e_2\beta=e_2\gamma=0$.
Hence we  have
\be
0=[e_2, \xi]\beta=(\nabla_{e_2}\xi-\nabla_{\xi}e_2)\beta=
(\gamma-\kappa_3)e_3\beta=3\gamma(\beta^2-1),\nonumber
\ee
which implies $\gamma=0$.

 The above argument shows that  the shape operator satisfies
 \be
 A\xi=\beta e_2, \ \ 
Ae_2=\beta\xi, \ \  Ae_3=0\nonumber 
 \ee 
at each point,  where $\beta\ne 0$. 
 Applying Lemma \ref{lem2}, we conclude that 
 $M$ is a minimal ruled real hypersurface.

{\bf Case (II)}. \ $n>2$.

Let $X$ be an arbitrary  vector field on $M$ perpendicular to $\xi$, $e_2$ and $e_3$. 
By Lemma \ref{lem1} we have $AX=\mu X$ and $A\phi X=\mu\phi X$.
Since  $\mu$ is constant, we obtain
\begin{align}
(\nabla_XA)\phi X-(\nabla_{\phi X} A)X &=\nabla_X(A\phi X)-A(\nabla_X(\phi X))-\nabla_{\phi X}(AX)+A(\nabla_{\phi X}X) \nonumber \\
&=(\mu-A)(\nabla_X(\phi X)-\nabla_{\phi X}X). \label{1}
\end{align}
Using (\ref{PA}) and (\ref{A}), we obtain
\begin{align}
g(\nabla_X(\phi X)-\nabla_{\phi X}X, \xi)&=-g(\nabla_X\xi, \phi X)+g(\nabla_{\phi X}\xi, X)\nonumber\\
&=-g(\phi AX, \phi X)+g(\phi A\phi X, X)\nonumber\\
&=-g(\phi AX, \phi X)-g(A\phi X, \phi X)\nonumber\\
&=-2\mu.\label{2}
\end{align}
By the equation (\ref{co}) of Codazzi for $Y=\phi X$, (\ref{A2}), (\ref{1}) and (\ref{2}), we get
\begin{align}
0&=g((\nabla_XA)\phi X-(\nabla_{\phi X} A)X, e_2)\nonumber \\
&=g(\nabla_X(\phi X)-\nabla_{\phi X}X, (\mu-A)e_2)\nonumber \\
&=(\mu-\gamma)g(\nabla_X(\phi X)-\nabla_{\phi X}X, e_2)-\beta
g(\nabla_X(\phi X)-\nabla_{\phi X}X, \xi)\nonumber \\
&=(\mu-\gamma)g(\nabla_X(\phi X)-\nabla_{\phi X}X, e_2)+2\beta\mu.\label{3}
\end{align}
Simliarly, we have
\begin{align}
-2c&=g((\nabla_XA)\phi X-(\nabla_{\phi X} A)X, \xi)\nonumber\\
&=g(\nabla_X(\phi X)-\nabla_{\phi X}X, (\mu-A)\xi)\nonumber\\
&=\gamma g(\nabla_X(\phi X)-\nabla_{\phi X}X, \xi)-\beta g(\nabla_X(\phi X)-\nabla_{\phi X}X, e_2)\nonumber\\
&=-2\mu\gamma-\beta g(\nabla_X(\phi X)-\nabla_{\phi X}X, e_2).\label{4}
\end{align}
Eliminating $g(\nabla_X(\phi X)-\nabla_{\phi X}X, e_2)$ from (\ref{3}) and (\ref{4}) gives 
\be 
(c-\mu\gamma)(\mu-\gamma)+\mu\beta^2=0. \label{Eq1}
\ee
Differentiating (\ref{Eq1}) along $e_3$ implies
\be
(2\mu\gamma-\mu^2-c)(e_3\gamma)+2\mu\beta(e_3\beta)=0. \label{Eq2}
\ee
Substituting (\ref{b}) and (\ref{r}) into (\ref{Eq2}), we have
\be
(2\mu\gamma-\mu^2-c)\{(2\gamma+\mu)\beta^2+2\gamma^3-\mu\gamma^2-(3\mu^2+c)\gamma+(2\mu^2+c)\mu\}=0.
\label{Eq3}
\ee
Computing the resultant of the left-hand sides of (\ref{Eq1}) and (\ref{Eq3}) with respect to 
$\beta$, we obtain
\be
(\gamma-\mu)(c\gamma-\mu^3)(2\mu\gamma-\mu^2-c)=0.
\label{Eq4}
\ee

If $\mu\ne 0$, then (\ref{Eq1}) and (\ref{Eq4}) show that
  $\beta$ and $\gamma$ are constant, that is, $M$ has constant principal curvatures.
  Let $h(p)$ the number of nontrivial projections of  $\xi_p$ onto the principal curvature spaces of $M $. Then, by Lemma \ref{lem1} we have $h=2$. According to the study of principal curvatures in \cite{diaz}, 
  $c=-1$ and $\alpha+\gamma=3\mu$ must be satisfied. However, this contradicts 
  algebraic condition in $(\ref{A})$.
  Hence, we obtain $\mu=0$, which implies $\gamma=0$ from
  (\ref{Eq4}). Therefore, by  Lemma \ref{lem2} and (\ref{A2}), we conclude that 
 $M$ is a minimal ruled real hypersurface.

Conversely,  it follows from Lemma \ref{lem1} and Lemma \ref{lem2} that a minimal ruled real hypersurface in
 $\tilde M^n(4c)$ attains  equality in $(\ref{ineq})$ at each point. 
The proof is finished. \qed

 \end{document}